\documentclass[11pt]{amsart}

\usepackage[colorlinks=true, pdfstartview=FitV, linkcolor=blue, 
citecolor=blue]{hyperref}

\usepackage{amssymb,amsmath,amscd}
\usepackage{bbm}
\usepackage{graphicx}
\usepackage{a4wide}
\newtheorem{theorem}{Theorem}[section]
\newtheorem{prop}[theorem]{Proposition}

\theoremstyle{definition}

\newcommand{\ts}{\hspace{0.5pt}}
\newcommand{\nts}{\hspace{-0.5pt}}

\newcommand{\RR}{\mathbb{R}\ts}
\newcommand{\CC}{\mathbb{C}}
\newcommand{\ZZ}{\mathbb{Z}}

\newcommand{\NN}{\mathbb{N}}

\newcommand{\QQ}{\mathbb{Q}}
\newcommand{\YY}{\mathbb{Y}}

\newcommand{\cF}{\mathcal{F}}
\newcommand{\cK}{\mathcal{K}}
\newcommand{\cL}{\mathcal{L}}

\newcommand{\vL}{\varLambda}

\newcommand{\dd}{\,\mathrm{d}}
\newcommand{\ee}{\ts\mathrm{e}}
\newcommand{\ii}{\mathrm{i}\ts}

\newcommand{\myfrac}[2]{\frac{\raisebox{-2pt}{$#1$}}
      {\raisebox{0.5pt}{$#2$}}}

\newcommand{\twovec}[2]{\bigl(\begin{smallmatrix} #1 \\ #2
      \end{smallmatrix}\bigr)}
\numberwithin{equation}{section}

\makeatletter
\renewcommand{\@captionfont}{\small}
\makeatother

\DeclareMathOperator{\dotcup}{\dot\cup}

\DeclareMathOperator{\dens}{dens}
\DeclareMathOperator{\sinc}{sinc}
\DeclareMathOperator{\vol}{vol}

\DeclareFontFamily{U}{mathx}{\hyphenchar\font45}
\DeclareFontShape{U}{mathx}{m}{n}{ <5> <6> <7> <8> <9> <10>
   <10.95> <12> <14.4> <17.28> <20.74> <24.88> mathx10 }{}
\DeclareSymbolFont{mathx}{U}{mathx}{m}{n}
\DeclareFontSubstitution{U}{mathx}{m}{n}
\DeclareMathAccent{\widecheck}{0}{mathx}{"71}

\newcommand{\sB}{\ts\underline{\nts B\!}\,}
\newcommand{\defeq}{\mathrel{\mathop:}=}
\newcommand{\eqdef}{=\mathrel{\mathop:}}

\begin{document}

\title[Variations on a theme by Fibonacci]
{Three variations on a theme by Fibonacci}

\author{Michael Baake}
\address{Fakult\"at f\"ur Mathematik, Universit\"at Bielefeld,
  \newline \indent Postfach 100131, 33501 Bielefeld, Germany}
\email{mbaake@math.uni-bielefeld.de}

\author{Natalie Priebe Frank}
\address{Department of Mathematics and Statistics, Vassar
  College,\newline \indent Poughkeepsie, NY 12604, USA}
\email{nafrank@vassar.edu}

\author{Uwe Grimm}
\address{School of Mathematics and Statistics, The Open
  University,\newline \indent Walton Hall, Milton Keynes MK7 6AA, UK}
\email{uwe.grimm@open.ac.uk}

\begin{abstract} 
  Several variants of the classic Fibonacci inflation tiling are
  considered in an illustrative fashion, in one and in two dimensions,
  with an eye on changes or robustness of diffraction and dynamical
  spectra. In one dimension, we consider extension mechanisms of
  deterministic and of stochastic nature, while we look at direct
  product variations in a planar extension. For the pure point part,
  we systematically employ a cocycle approach that is based on the
  underlying renormalisation structure. It allows explicit
  calculations, particularly in cases where one meets regular model
  sets with Rauzy fractals as windows.
\end{abstract}

%\keywords{inflation tilings, dynamical spectrum, diffraction, cocycle}
%\subjclass{42B10, 52C23, 37F25}

\maketitle

\centerline{Dedicated to Manfred Denker on the occasion of his 75th
  birthday}

\section{Introduction}

The mathematical theory of aperiodic order profits enormously from two
construction principles, namely the cut and project method and the
substitution or inflation method; see \cite{TAO,Nat-primer} and
references therein for general background. While the understanding of
cut and project sets has reached a rather satisfactory level, this is
less so for primitive inflation tilings. In particular, the
classification of their spectra is as yet incomplete. While the still
open Pisot substitution conjecture \cite{Aki,Bernd} marks the frontier
for one-dimensional systems, the situation is worse in two or more
dimensions. Here, it is already less clear what an appropriate
conjecture could be, because there are rather intricate additional
constraints of a more \emph{geometric} origin that do not exist in one
dimension.

In this exposition, we reconsider the simplest and best-studied
example, the Fibonacci chain, and investigate modifications as well as
extensions to planar systems. On the line, we demonstrate two
mechanisms that add continuous spectral components due to disorder,
either of random or of deterministic type.  In the case of planar
tilings, we are particularly interested in geometric variations that
change the system topologically, but not measure-theoretically, which
is to say that we are probing the stability of spectral properties.

Our starting point is the self-similar Fibonacci \emph{tiling
  dynamical system} (TDS) in one dimension, as defined by the
primitive inflation rule
\[
    \varrho : \; a \mapsto ab \, , \; b \mapsto a \ts ,
\]
with $a$ and $b$ considered as \emph{tiles} (or intervals) of length
$\tau = \frac{1}{2} \bigl( 1 + \mbox{\small $\sqrt{5}$} \, \bigr)$ and
$1$, respectively. The substitution matrix of $\varrho$ is
\begin{equation}\label{eq:fibmat}
    M \, = \, \begin{pmatrix} 1 & 1 \\ 1 & 0 \end{pmatrix},
\end{equation}
which is primitive, with Perron--Frobenius (PF) eigenvalue $\tau$ and
corresponding left and right eigenvectors, in Dirac notation,
\[
  \langle u\ts | \, = \, \myfrac{\tau+2}{5} \bigl(\tau,1\bigr)\quad
\text{and} \quad
   |\ts v\rangle \, = \, \bigl(\tau^{-1},\tau^{-2}\bigr)^{T}.
\]
They are normalised such that $\langle 1 | \ts v \rangle =1$, which
means that the entries of $|\ts v\rangle$ are the relative frequencies
of the tiles, together with $\langle u\ts |\ts v\rangle=1$. In this
setting, one has
\begin{equation}\label{eq:Pdef}
    \lim_{n\to\infty} \tau^{-n} M^n \, = \, \myfrac{\tau+2}{5}
    \begin{pmatrix} 1 & \tau^{-1} \\ \tau^{-1} & \tau^{-2}\end{pmatrix}
    \, = \, |\ts v\rangle\langle u\ts | \, \eqdef \, P\ts ,
\end{equation}
where $P=P^2$ is a symmetric projector of rank $1$ with 
spectrum $\{1,0\}$. 

If $\YY$ denotes the (compact) tiling hull defined by $\varrho$, the
topological dynamical system $(\YY, \RR)$ is strictly ergodic; compare
\cite[Ch.~5]{DGS}.  More precisely, $\YY$ can be constructed as the
orbit closure of a bi-infinite Fibonacci tiling, for which one usually
employs one of the two fixed points of $\varrho^2$, with core $a|a$ or
$b|a$; see \cite[Ex.~4.6]{TAO} for details and \cite{Weiss} for
general background. Since each fixed point is repetitive, the
resulting dynamical system is minimal.  The unique probability measure
on $\YY$ is the \emph{patch frequency measure}, denoted by $\mu$.  Via
the left endpoints of the intervals, one can consider the elements of
$\YY$ either as tilings or as Delone sets, which are two viewpoints
that we will tacitly identify. This is justified by the fact that the
two structures are \emph{mutually locally derivable} (MLD) from one
another (and thus certainly topologically conjugate); see
\cite[Sec.~5.2]{TAO} for background. Note that, in our setting, each
such Delone set has density $(\tau+2)/5=\tau/\sqrt{5}$.

The topological dynamical system $(\YY, \RR)$ is uniquely ergodic and
has pure point spectrum in the measure-theoretic sense, employing the
Koopman operator on the Hilbert space $L^2 (\YY,\mu)$. All
eigenfunctions have continuous representatives \cite[Thm.~1.4]{Host},
though the system itself is not equicontinuous. In the length scale
chosen, the spectrum is
$L^{\circledast}=\ZZ [\tau]/\sqrt{5}\subset\QQ(\sqrt{5}\,)$, in
additive notation, as can easily be extracted from the description of
$\YY$ by the projection method; see \cite[Sec.~9.4.1]{TAO}. The
understanding of this example is based on its simultaneous description
as an inflation tiling and as a cut and project set.  We will recall
the details in Section~\ref{sec:Fibo}, based on the Minkowski
embedding of $\ZZ [\tau]$ as a lattice in $\RR^2$.

Below, we go through a number of modifications or variations, starting
with one-dimensional examples with some added structure, random or
deterministic, to show origins of continuous spectral components of
different type.  Then, we double the dimension and consider
Fibonacci-type tilings of the plane, and some of their rearrangements
via modified inflation rules. Here, we present two rather typical
phenomena, namely the robustness of pure point spectra in the Pisot
case and the emergence of Rauzy fractals even in simple
situations.\smallskip

Our exposition has a somewhat informal character, aiming at an
illustration of several possible phenomena via simple examples. To do
so, we present the various derivations first, and then summarise our
results in formal theorems. We begin with a brief summary of the
Fibonacci TDS in Section~\ref{sec:Fibo}, followed by two variations in
Section~\ref{sec:first} that add absolutely continuous or singular
continuous components.  Then, in Section~\ref{sec:second}, we define
the natural direct product Fibonacci TDS in the plane, and derive its
properties.  Finally, we investigate \emph{direct product variations}
(DPVs) in Section~\ref{sec:third}, which show robustness of the
spectral type.

\section{The Fibonacci chain and its spectral 
properties}\label{sec:Fibo}

Here, we begin with a condensed summary of the standard derivation of
the spectrum, with some focus on the eigenfunctions. Then, we show how
to use the renormalisation approach to obtain the same results, which
will be employed for our variations.

\subsection{Standard approach}
The fact that $(\YY,\RR)$ has pure point (dynamical) spectrum is
equivalent to the statement that the Fibonacci chain has pure point
diffraction; see \cite{BL} and references therein for the general
theory. Here, for any Delone set $\vL\in\YY$, one defines the measure
$\omega=\delta^{}_{\! \vL} \defeq \sum_{x\in\vL} \delta^{}_{x}$,
called the \emph{Dirac comb} of $\vL$. Then, $\omega$ leads to the
autocorrelation $\gamma = \omega \circledast \widetilde{\omega}$,
where $\circledast$ denotes volume-averaged (or Eberlein) convolution
\cite[Sec.~8.8]{TAO}, and to the diffraction measure
$\widehat{\gamma}$, which is the (existing) Fourier transform of the
autocorrelation.

For the special case of the Fibonacci chain, one obtains the pure 
point measure
\[
    \widehat{\gamma} \, = \sum_{k\in L^{\circledast}} I(k) \, \delta^{}_{k}
    \quad \text{with} \quad
    I (k) \, = \, \lvert A^{}_{\!\vL} (k) \rvert^2 ,
\]
where the amplitudes, or \emph{Fourier--Bohr} (FB) coefficients, are
given by
\[
     A^{}_{\!\vL} (k) \, = \lim_{r\to\infty} \myfrac{1}{2 r} 
     \sum_{x\in\vL_{r}} \ee^{-2\pi\ii k x},
\]
with $\vL_{r} \defeq \vL \cap [-r,r]$. The limits exist uniformly, and
one has
\[
     A^{}_{t+\vL} (k) \, = \, \ee^{-2\pi\ii k t} A^{}_{\!\vL} (k)  \ts ,
\]
which means that, for any $k$ for which the coefficient is
non-trivial, $\vL \mapsto A^{}_{\!\vL} (k)$ defines an eigenfunction
of our system, with eigenvalue $k$ in additive notation. Since
eigenfunctions for primitive inflation tilings have continuous
representatives \cite{Host,Boris}, their knowledge on the defining
orbit suffices to determine them.

For later convenience, we consider the two points sets $\vL_{a}$ and
$\vL_{b}$ of left endpoints of tiles of types $a$ and $b$ separately,
with $\vL = \vL_{a} \ts\ts\dot{\cup}\ts\ts \vL_{b}$, and define
corresponding amplitudes $A^{}_{\!\vL_{a}}$ and $A^{}_{\!\vL_{b}}$,
where $A^{}_{\!\vL}=A^{}_{\!\vL_{a}}\! +A^{}_{\!\vL_{b}}$. When $\vL$
is one of the two fixed points mentioned above and $k\in
L^{\circledast}$, the amplitudes are given \cite[Sec.~9.4.1]{TAO} by
\begin{equation}\label{eq:1Dampli}
\begin{split}
   A^{}_{\!\vL_{a}}(k) \, &= \, 
   \frac{1}{\sqrt{5}} \int_{\tau-2}^{\tau-1} \ee^{2\pi\ii k^{\star} y} \dd y
   \, = \, \myfrac{\ee^{2\pi\ii k^{\star}(\tau-1)}-
     \ee^{2\pi\ii k^{\star}(\tau-2)}}
      {2\pi\ii \sqrt{5}\ts k^{\star}_{\phantom{I}}} \ts , \\[1.5mm]
   A^{}_{\!\vL_{b}}(k) \, &= \, 
   \frac{1}{\sqrt{5}} \int_{-1}^{\tau-2} \ee^{2\pi\ii k^{\star} y} \dd y
   \, = \, 
   \myfrac{\ee^{2\pi\ii k^{\star}(\tau-2)}-\ee^{-2\pi\ii k^{\star}}}
                {2\pi\ii \sqrt{5}\ts k^{\star}_{\phantom{I}}} \ts ,
\end{split}
\end{equation}
while they vanish for all other values of $k$. Here, $k^\star$ is the
image of $k$ under the $\star$-map which acts as algebraic conjugation
with $\sqrt{5}\mapsto -\sqrt{5}$ on $\QQ(\sqrt{5}\,)$.  The explicit
formulas for the amplitudes emerge from the description of the
Fibonacci chain as a regular model set. Here, we use the natural
\emph{cut and project scheme} (CPS) with the lattice
\begin{equation}\label{eq:Minkowski}
  \cL \, \defeq  \, \bigl\{ (x,x^{\star}) :
     x \in \ZZ [ \tau] \bigr\} \ts ,
\end{equation}   
which is the standard Minkowski embedding of $\ZZ [\tau]$ in $\RR^2$;
see \cite[Fig.~3.3]{TAO} for an illustration and \cite[Sec.~7.2]{TAO}
for further details on the projection formalism in this setting.

With this approach, one finds\footnote{Here, we use $\widehat{g}$ for
  the Fourier transform of a function $g$, and $\widecheck{g}$ for its
  inverse transform.}
\begin{equation}\label{eq:Fampli}
    A^{}_{\!\vL_{a,b}}(k)\, = \, \frac{\dens(\vL_{a,b})}{\vol(W_{\! a,b})}\,
    \widehat{1^{}_{W_{\! a,b}}}(-k^{\star}) \, = \, 
    \myfrac{1}{\sqrt{5}} \,\widecheck{1^{}_{W_{\! a,b}}}(k^{\star})\ts ,
\end{equation}
where 
\begin{equation}\label{eq:Fib-windows}
   W_{\! a} \, = \, [\tau-2,\tau-1] 
   \quad \text{and} \quad
    W_{b} \, = \, [-1,\tau-2]
\end{equation}   
are the closures of the windows for the point sets $\vL_{a}$ and
$\vL_{b}$ in the projection formalism. In fact, one has the true
inclusion
\[
     \vL^{}_{a,b} \, \subset \, \{ x\in \ZZ [\tau] :
        x^{\star} \in W^{}_{a,b} \} \ts ,  
\]
where the sets on the right-hand side contain one extra point each,
caused by one of the window boundary points; see \cite[Ex.~7.3]{TAO}
for details. This fine point is important for the topological
structure of the hull, but of no relevance to the spectral
considerations.

For the sum of the amplitudes, the formulas simplify to
\[
     A^{}_{\!\vL} (k) \, = \, \begin{cases} \frac{\tau}{\sqrt{5}}\ts
        \ee^{\pi \ii k^{\star} (\tau-2)} \sinc(\pi\tau k^{\star}), 
            & \text{if $k\in L^{\circledast}$}, \\
        0 , & \text{otherwise}, \end{cases}
\]
where $\sinc(x)=\sin(x)/x$. Consequently, the diffraction intensities
are
\[
    I(k) \, = \,  \lvert A^{}_{\!\vL}(k)\rvert^{2}\, = \, \begin{cases}
    \Bigl( \frac{\tau}{\sqrt{5}} \sinc (\pi \tau k^{\star}) \Bigr)^{2},
            & \text{if $k\in L^{\circledast}$}, \\
        0 , & \text{otherwise}. \end{cases}
\]
Note that the intensity function, and thus the diffraction measure
$\widehat{\gamma}$, is the same for all $\vL\in\YY$, while the
amplitudes will usually differ by a phase factor for distinct
elements.

More generally, if we introduce two (in general complex) weights for
the two different endpoints ($u^{}_{a}$ and $u^{}_{b}$ say, so
$\omega = u^{}_{a}\delta^{}_{\!\vL_{a}}+u^{}_{b}\delta^{}_{\!\vL_{b}}$),
we obtain the intensity by superposition as
\begin{equation}\label{eq:Fibo-int}
  I(k) \, = \, \widehat{\,\gamma^{}_{(u^{}_{a},u^{}_{b})}}
  \bigl(\{k\}\bigr)\, = \, 
    \bigl|  u^{}_{a} \, A^{}_{\!\vL_{a}}(k) + 
               u^{}_{b} \, A^{}_{\!\vL_{b}}(k)\bigr|^{2}.
\end{equation}
Let us sum up these well-known spectral properties \cite{Q,TAO} as
follows.

\begin{theorem}
  The Fibonacci dynamical system\/ $(\YY,\RR)$, in its geometric
  realisation as described above, is strictly ergodic and has pure
  point spectrum, both in the diffraction and in the dynamical
  sense. For given weights\/ $u^{}_{a}, u^{}_{b} \in \CC$, the
  diffraction measure is given by
\[
  \widehat{\,\gamma^{}_{(u^{}_{a},u^{}_{b})}} \, = \sum_{k\in
    L^{\circledast}} I(k) \, \delta^{}_{k}
\]
with the intensities\/ $I (k)$ from Eq.~\eqref{eq:Fibo-int}. The
Fourier module is\/ $L^{\circledast} = \ZZ [\tau]/\sqrt{5}$ and agrees
with the dynamical spectrum of\/ $(\YY, \RR)$ in additive notation.
\qed
\end{theorem}
             
Note that the autocorrelation measure
$\gamma^{}_{(u^{}_{a},u^{}_{b})}$, which is a pure point measure with
Meyer set\footnote{A point set $\vL \subset \RR$ is a \emph{Meyer set}
  if it is relatively dense and if $\vL - \vL$ is uniformly discrete.}
support, can be expressed in terms of the (dimensionless) pair
correlation coefficients
\[
    \nu_{\alpha\beta}(z) \, \defeq \, 
    \frac{\dens\bigl(\vL_{\alpha}\cap(\vL_{\beta}-z)\bigr)}
    {\dens(\vL)} \, = \, \nu^{}_{\beta\alpha} (-z) \ts ,
\]
which are positive for all $z\in\vL_{\beta}-\vL_{\alpha}$ and $0$
otherwise, as
\[
    \gamma^{}_{(u^{}_{a},u^{}_{b})}\bigl(\{z\}\bigr) \, = \, 
    \dens(\vL) \sum_{\alpha,\beta\in\{a,b\}} \overline{u^{}_{\alpha}}
    \: \nu^{}_{\alpha\beta}(z)\, u^{}_{\beta}.
\]
In particular, we have $\nu^{}_{aa} (0) = \tau^{-1}$ and
$\nu^{}_{bb} (0) = \tau^{-2}$, hence $\nu^{}_{aa} (0)
+ \nu^{}_{bb} (0) = 1$.

\subsection{Renormalisation-based approach}
While the relation between the FB coefficients and the Fourier
transform of the compact windows holds for regular model sets in
general, see \cite[Thm.~9.4]{TAO}, it is practically impossible to
compute the coefficients by Fourier transform of the windows if the
latter are compact sets with fractal boundaries. Let us therefore
explain a different approach that will also work in such more
complicated situations.

With $\sigma \defeq \tau^{\star}=1-\tau$, the inflation structure
induces a relation between the windows that, in terms of their
characteristic functions, reads
\[
    1^{}_{W_{\! a}} \, =\,  1^{}_{\sigma W_{\! a}\cup\ts\sigma W_{b}} 
    \quad\text{and}\quad
    1^{}_{W_{b}} \, = \, 1^{}_{\sigma W_{\! a} +\ts\sigma} \ts ,
\]
as can easily be verified for the windows from
\eqref{eq:Fib-windows}. Observing that
$1^{}_{\sigma W_{\! a}\cup\ts\sigma W_{b}} =1^{}_{\sigma W_{\! a}} +
1^{}_{\sigma W_{b}}$ holds as an equation of $L^{1}$-functions, an
application of the inverse Fourier transform gives
\begin{equation}\label{eq:charfun}
   \widecheck{1^{}_{W_{\! a}}} \, = \, 
   \widecheck{1^{}_{\sigma W_{\! a}}} + \widecheck{1^{}_{\sigma W_{b}}}
    \quad\text{and}\quad
   \widecheck{1^{}_{W_{b}}} \, = \, 
   \widecheck{1^{}_{\sigma W_{\! a}+\ts \sigma}}\ts .
\end{equation}
By an elementary calculation, one finds
\begin{equation}\label{eq:affineFT}
    \widecheck{1^{}_{\alpha K+\beta}}(y) \, = \, 
    \lvert\alpha\rvert \, \ee^{2\pi\ii \beta y} \, 
    \widecheck{1^{}_{K}}(\alpha\ts y)\ts ,
\end{equation}
which holds for arbitrary $\alpha,\beta\in\RR$ with $\alpha\ne 0$ 
and any compact set $K\subset \RR$. 

Defining $h^{}_{a,b}=\widecheck{1^{}_{W_{\! a,b}}}$, an application of
\eqref{eq:affineFT} to Eq.~\eqref{eq:charfun} results in
\begin{equation}\label{eq:self}
   \begin{pmatrix} h^{}_{a} \\ h^{}_{b} \end{pmatrix} (y) \, = \, 
   \lvert\sigma\rvert  \, \sB (y) 
   \begin{pmatrix} h^{}_{a} \\ h^{}_{b} \end{pmatrix} (\sigma y)
\quad\text{with}\quad
\sB (y) \, \defeq \, \begin{pmatrix} 1 & 1 \\ 
    \ee^{2\pi\ii\sigma y} & 0 \end{pmatrix},
\end{equation}
where $\sB$ is related to the Fourier matrix from the renormalisation
approach \cite{BG,BFGR} by first taking the $\star$-map of the
set-valued displacement matrix
$T=\left(\begin{smallmatrix} \{0\} & \{0\} \\ \{\tau\} & \varnothing
\end{smallmatrix}\right)$ and then its (inverse) Fourier
transform. For this reason, we call it the \emph{internal Fourier
  matrix} \cite{BG19}. The latter can also be viewed as emerging from
the commutative diagram
\begin{equation}\label{eq:diagram}
\begin{CD}
T @>\cF^{-1}>> B(.)\\
@V{\star}VV @VV{\mbox{\tiny\textcircled{$\star$}}}V\\
T^{\ts\star} @>>\cF^{-1}> \sB(.)
\end{CD}
\end{equation}
where $\cF$ denotes the Fourier transform of a matrix of (finite)
Dirac combs and \raisebox{2pt}{\tiny\textcircled{$\star$}} the induced
mapping on the level of the Fourier matrices.

Using the notation $|\ts h\rangle=(h^{}_{a},h^{}_{b})^T$ and applying
the above iteration $n$ times leads to
\[
    |\ts h(y)\rangle \, = \, \lvert\sigma\rvert^{n} \sB^{(n)}(y) \,
    |\ts h(\sigma^{n} y)\rangle \quad \text{where}\quad
    \sB^{(n)}(y)\, \defeq \, \sB (y) \sB (\sigma y) \cdots 
    \sB (\sigma^{n-1} y)\ts ,
\]
which can be interpreted as a singular case of a transfer matrix
approach to a cocycle; compare \cite{AD,Pohl}.  In particular,
$\sB^{(1)}=\sB$ and $\sB^{(n)}(0)=M^{n}$ for all $n\in\NN$, where $M$
is the substitution matrix from \eqref{eq:fibmat}. Note that
$\sB^{(n)}(y)$ defines a matrix cocycle, called the \emph{internal
  cocycle}, which by \eqref{eq:diagram} is related to the usual
inflation cocycle by an application of the $\star$-map to the
displacement matrices of the powers of the inflation rule
\cite{BGM,BG19}.

It is not difficult to prove that $|\ts h(y)\rangle=C(y) |\ts
h(0)\rangle$, where
\[
   C(y) \, \defeq \, \lim_{n\to\infty} \lvert\sigma\rvert^n \sB^{(n)}(y)
\]
exists pointwise for every $y\in\RR$. In fact, one has compact
convergence, which implies that $C(y)$ is continuous \cite[Thm.~4.6
and Cor.~4.7]{BG19}.  For any $m,n\in\NN$, one has
\begin{equation}\label{eq:Bsplit}
   \sB^{(n+m)}(y) \, = \, \sB^{(n)}(y) \,\sB^{(m)}(\sigma^n y)\ts .
\end{equation}
Employing this relation with $m=1$, letting $n\to\infty$, and
observing $\lvert\sigma\rvert=\tau^{-1}$, one obtains
\[
   \tau \, C(y) \, = \, C(y) M\ts .
\]
This relation implies that each row of $C(y)$ is a multiple of
$\langle u|$, so we may define a vector-valued function
$|\ts c(y) \rangle$ such that $C(y)=|\ts c(y)\rangle\langle u\ts |$
holds with $|\ts c(y)\rangle = \bigl(c^{}_{a}(y),c^{}_{b}(y)\bigr)^T$.
Since $C(0)=P$ with the projector $P=|\ts v\rangle\langle u\ts |$ from
\eqref{eq:Pdef}, we have $|\ts c(0)\rangle=|\ts v\rangle$.

As $|\ts h(y)\rangle=|\ts c(y)\rangle\langle u\ts |\ts h(0)\rangle$,
where $|\ts h(0)\rangle=\tau \ts |\ts v\rangle$ follows from a simple
calculation, we get
\[
    |\ts h(y)\rangle \, = \, \tau \ts\ts |\ts c(y)\rangle\ts ,
\]
and the inverse Fourier transforms of the windows are encoded in the
matrix $C$. For the Fibonacci case at hand, we can explicitly
calculate $|\ts c(y)\rangle$ through the Fourier transforms of the
known windows to be
\[
     c^{}_{a}(y) \, = \, 
     \myfrac{\ee^{2\pi\ii (\tau-1)y} - \ee^{2\pi\ii (\tau-2)y}}
     {2\pi\ii y}   \quad\text{and}\quad
     c^{}_{b}(y) \, = \, 
     \myfrac{\ee^{2\pi\ii (\tau-2)y} - \ee^{-2\pi\ii y}}
     {2\pi\ii y} \ts ,
\]
from which one can explicitly check that there is no $y\in\RR$ for
which both functions vanish simultaneously. Consequently, $C(y)$ is
always a matrix of rank $1$.

In other situations, where no explicit formula for the Fourier
transforms is available, one can replace $C(y)$ by
$\lvert \sigma\rvert^n \sB^{(n)}(y)$ or its analogue for a
sufficiently large $n$, subject to the condition that $C(y)$ is
approximated sufficiently well (in some matrix norm, say) and that
$\lvert \sigma\rvert^n y$ is close enough to $0$. This works because
the windows are compact sets, so that their Fourier transforms are
continuous functions. The convergence of this approximation turns out
to be exponentially fast; see \cite{BG19} for details and an extension
of the cocycle approach to more general inflation systems.

With our previous relation \eqref{eq:Fampli}, for $k\in
L^{\circledast}$, the FB amplitudes are
\[
   A^{}_{\!\vL_{ a,b}}(k) \, = \, 
   \frac{h^{}_{a,b}(k^{\star})}{\sqrt{5}} \ts , 
\]
which means that they can now be calculated via $C$ as well. This
gives us a way to compute the eigenfunctions and the general
diffraction amplitudes to arbitrary precision.  Unless stated
otherwise, numerical calculations and illustrations below will always
be based on the cocycle approach due to its superior speed and
accuracy in the presence of complex windows.

\section{Two variations: Randomness and deterministic
  disorder}\label{sec:first}

The goal of this section is to demonstrate two mechanisms that can
alter the spectrum by adding a continuous component.

\subsection{Randomness} Given a Fibonacci tiling from the hull $\YY$,
we now introduce some uncorrelated disorder into the $a$ positions,
by randomly assigning two different labels or weights. This way, we
address and answer a question by Strungaru \cite{Nicu-PC} on the
spectral consequences of this type of modification.  Consider the
situation where tiles of type $a$ are independently replaced by a tile
$\underline{a}$ of the same length with probability $q$, and kept with
probability $p=1-q$. By a simple calculation, whenever $p\in (0,1)$,
this leads to a shift space with topological entropy $\log (2) /\tau$,
calculated per tile (rather than per unit length).

By the strong law of large numbers, using the results of \cite{BBM},
the pair correlation coefficients for this modified system can be
expressed in terms of the original coefficients as follows,
\begin{align*}
    \nu^{(p)}_{a\underline{a}}(z) \, =\,
      \nu^{(p)}_{\underline{a}a}(z) &\,= \,
   \begin{cases} p \ts q\, \nu^{}_{aa}(z), & z\ne 0,\\
          0, & z=0,\end{cases} \qquad \\[1mm]
\nu^{(p)}_{aa}(z) \,  = \,
     \begin{cases} p^2 \nu^{}_{aa}(z), & z\ne 0,\\
        p/\tau, & z=0,\end{cases}\; & \qquad
        \nu^{(p)}_{\underline{a}\underline{a}}(z) \,= \,
     \begin{cases} q^2 \nu^{}_{aa}(z), & z\ne 0,\\
        q/\tau, & z=0,\end{cases}\\[1mm]
  \nu^{(p)}_{ab}(z)\, =\, \nu^{(p)}_{ba}(-z)
       \, = \,  p\, \nu^{}_{ab}(z), \! & \qquad
  \nu^{(p)}_{\underline{a}b}(z)\, =\, \nu^{(p)}_{b\underline{a}}(-z) 
    \, = \,  q\, \nu^{}_{ab}(z),\\[2mm]
   \nu^{(p)}_{bb}(z) &\, = \,  \nu^{}_{bb}(z).
\end{align*}
Now, we have three different types of points. Using weights
$u^{}_{a}$, $u^{}_{\underline{a}}$ and $u^{}_{b}$, we thus obtain the
autocorrelation
\[
   \gamma^{(p)}_{(u^{}_{a},u^{}_{\underline{a}},u^{}_{b})}
     \bigl(\{z\}\bigr) \, = \, 
     \dens(\vL) \!\sum_{\alpha,\beta\in\{a,\underline{a},b\}}\!
    \overline{u^{}_{\alpha}}\: \nu^{(p)}_{\alpha\beta} (z)
    \, u^{}_{\beta}.
\]
Setting $v^{}_{a}=p\ts u^{}_{a}+qu^{}_{\underline{a}}$ and
$v^{}_{b}=u^{}_{b}$, an explicit computation shows that,
for $z\ne 0$, this becomes 
\[
    \gamma^{(p)}_{(u^{}_{a},u^{}_{\underline{a}},u^{}_{b})}
   \bigl(\{z\}\bigr) \, = \, 
   \dens(\vL) \sum_{\alpha,\beta\in\{a,b\}}
    \overline{v^{}_{\alpha}}\: \nu^{}_{\alpha\beta} (z)
    \, v^{}_{\beta} \, = \,
    \gamma^{}_{(v^{}_{a},v^{}_{b})}\bigl(\{z\}\bigr),
\]
while, for $z=0$, we obtain
\begin{align*}
      \gamma^{(p)}_{(u^{}_{a},u^{}_{\underline{a}},u^{}_{b})}
     \bigl(\{0\}\bigr) \, & = \, 
      \dens(\vL) \Bigl(\frac{p}{\tau} \ts |u^{}_{a}|^2
      + \frac{q}{\tau} \ts |u^{}_{\underline{a}}|^2 +
      \myfrac{1}{\tau^2} \ts |u^{}_{b}|^2\Bigr)\\[2mm]
      & = \, \gamma^{}_{(v^{}_{a},v^{}_{b})}\bigl(\{0\}\bigr)
        \, +\, p\ts q \ts \ts \bigl|u^{}_{a}-u^{}_{\underline{a}}
        \bigr|^2 \dens(\vL) \ts .
\end{align*}
Taking the Fourier transform, we see that the additional point measure
at $z=0$ gives rise to an absolutely continuous component, and we
obtain the following result, which can also be seen as a special case
of the random cluster model treated in \cite{BBM}.

\begin{theorem}
  If the Fibonacci chain is randomised at the\/ $a$-positions with a
  binary Bernoulli process with probabilities\/ $p\in [0,1]$ and\/
  $q=1-p$ as described above, the diffraction measure is almost surely
  given by
\[  
  \widehat{\gamma^{(p)}_{(u^{}_{a},u^{}_{\underline{a}},u^{}_{b})}
    \hspace{-2em}}
    \hspace{2em} \, = \, 
   \widehat{\gamma^{\phantom{\scriptscriptstyle p}}_{(v^{}_{a}, v^{}_{b})}}
   \, + \, p \ts q \, \bigl|u^{}_{a}-u^{}_{\underline{a}}\bigr|^2
   \lambda^{}_{\mathrm{L}} \ts ,
\]  
with\/ $v^{}_{a}=pu^{}_{a}+qu^{}_{\underline{a}}$ and\/
$v^{}_{b}=u^{}_{b}$, where\/ $\lambda^{}_{\mathrm{L}}$ denotes
Lebesgue measure on $\RR$. In particular, for\/ $p\ts q \ne 0$ and\/
$u^{}_{a} \ne u^{}_{\underline{a}}$, it is a sum of a pure point and
an absolutely continuous measure. \qed
\end{theorem}

Notice that the absolutely continuous component vanishes if
$u^{}_{a}=u^{}_{\underline{a}}$ or $p\ts q = 0$. Then, we recover the
result for the perfect Fibonacci chain in this case, for two different
reasons. When $p\ts q=0$, almost surely only one type of $a$ is
present, while $u^{}_{a} = u^{}_{\underline{a}}$ make the two types
indistinguishable from a diffraction point of view.

\subsection{Deterministic disorder}
Inspired by \cite[Sec.~8]{KS} and \cite{BG19}, let us move on to our
second variation and extend the Fibonacci substitution to a
four-letter substitution with bar-swap symmetry \cite{BG}, namely
\begin{equation}\label{eq:twist-Fib}
    a\mapsto ab\ts , \quad 
    \underline{a}\mapsto \underline{a}\underline{b}\ts ,\quad
    b\mapsto \underline{a}\ts ,\quad
    \underline{b}\mapsto a\ts .
\end{equation}
The substitution matrix is
\[
    M \, = \, \begin{pmatrix} 1 & 0 & 0 & 1 \\
    0 & 1 & 1 & 0 \\ 1 & 0 & 0 & 0 \\ 0 & 1 & 0 & 0 
    \end{pmatrix}
\]
with spectrum
$\bigl\{\tau,1-\tau, \frac{1}{2}\bigl(1\pm \ii
\sqrt{3}\,\bigr)\bigr\}$ and the PF eigenvectors
$\langle u\ts | = \frac{2(3\tau+4)}{19} (\tau,\tau,1,1)$ and
$|\ts v\rangle = \frac{2-\tau}{2} \ts (\tau,\tau,1,1
)^T$. Normalisation is again such that
$\langle u\ts|\ts v\rangle = \langle 1| \ts v \rangle = 1$.

Let us choose natural interval lengths $\tau$ for $a, \underline{a}$
and $1$ for $b,\underline{b}$, which means that we use the same
setting as for the perfect Fibonacci chain, including the ring of
integers, $\ZZ[\tau]$, and its Minkowski embedding $\cL$ from
\eqref{eq:Minkowski}. The one-sided fixed point equations for the
Delone sets of left endpoints read
\[
     \vL_{a} \, = \, \tau\vL_a \dotcup \tau \vL_{\underline{b}}\ts, \quad
     \vL_{\underline{a}} \, = \, \tau\vL_{\underline{a}} 
     \dotcup \tau \vL_{b}\ts, \quad
     \vL_{b} \, = \, \tau\vL_{a}+\tau\ts,\quad
      \vL_{\underline{b}} \, = \, \tau\vL_{\underline{a}}+\tau\ts .
\]
Applying the $\star$-map and taking closures leads to
\[
     W_{\! a} \, = \, \sigma W_{\! a} \cup 
     \sigma  W_{\nts\underline{b}}\ts, \quad
     W_{\!\underline{a}} \, = \, \sigma W_{\!\underline{a}} 
     \cup \sigma  W_{\nts b}\ts, \quad
     W_{\nts b} \, = \, \sigma W_{\!a}+\sigma \ts,\quad
      W_{\nts\underline{b}} \, = \, \sigma W_{\!\underline{a}}+\sigma \ts .
\]
This defines a contractive IFS on $(\cK \RR)^4$, where
$\mathcal{K} \RR$ is the space of non-empty compact subsets of $\RR$.
This is a complete metric space when equipped with the Hausdorff
metric as distance between compact sets. By the contraction principle,
also known as Hutchinson's theorem in this context, the IFS has a
unique solution \cite[Thm.~1.1]{BM-self}; see \cite{DY} for related
results.

One can verify that this solution is given by
\[
   W_{\! a} \, = \, W_{\!\underline{a}} \, = \, [\tau-2,\tau-1]\quad
   \text{and}\quad
   W_{\nts b} \, = \, W_{\nts\underline{b}} \, = \, [-1,\tau-2]\ts .
\]
The (perhaps surprising) point here is that these windows define true
covering sets of the four types of Delone sets that are immensely
useful for the spectral analysis. In particular, as detailed in
\cite{BG19}, one can extract the pure point part of the spectrum from
them.

Since these windows are the ones of the original Fibonacci chain from
\eqref{eq:Fib-windows}, we see that the disjoint point sets
$\vL_{a}^{\star}$ and $\vL_{\underline{a}}^{\star}$, as well as
$\vL_{b}^{\star}$ and $\vL_{\underline{b}}^{\star}$, are dense in the
same sets. In particular, this shows that none of the individual point
sets $\vL_{\alpha}$, with
$\alpha\in\{a,\underline{a},b,\underline{b}\}$, is a model set. Note
that the bar-swap inflation \eqref{eq:twist-Fib} leads to a
deterministic but non-trivial partition of the $a$ and $b$ positions
of the original Fibonacci chain into two sets each.

What is more, the images of the four types of point sets under the
$\star$-map are still uniformly distributed in the corresponding
windows \cite[Thm.~5.3]{BG19}. This implies that the FB coefficients,
by the standard Weyl equidistribution argument \cite[Lemma~9.4]{TAO},
are given by
\[
   A_{\alpha}(k) \, = \, \myfrac{1}{2\sqrt{5}}
   \, \widecheck{1^{}_{W_{\!\alpha}}}(k^{\star})\ts .
\]
Since we know that the Bombieri--Taylor property holds for primitive
inflation systems, see \cite[Thm.~3.23 and Rem.~3.24]{BGM}, the pure
point part of the diffraction of a weighted Dirac comb
$\sum_{\alpha}u^{}_{\alpha}\,\delta^{}_{\!\vL_{\alpha}}$ is a sum of
the form $\sum_{k\in L^{\circledast}}I(k)\ts\delta^{}_{k}$ with
\begin{equation}\label{eq:twist-intens}
  I(k) \, = \, \Bigl|\textstyle{\sum\limits_{\alpha}}
  \, u^{}_{\alpha}\ts A^{}_{\alpha}(k)\Bigr|^{2}.
\end{equation}

As the original Fibonacci system is a factor of its twisted bar-swap
extension, via identifying $a$ with $\underline{a}$ and $b$ with
$\underline{b}$, it is clear that the dynamical spectrum of our
bar-swap extension is of mixed type, which will be reflected for the
generic choice of the weights $u^{}_{\alpha}$ in the diffraction
measure as well. Since the extension is (measure-theoretically)
\mbox{2{\ts\ts}:1}, as follows from ordinary Fibonacci being an almost
everywhere \mbox{2{\ts\ts}:1} factor, we know that the continuous part
of the spectrum must be of pure type \cite{BG}, which turns out to be
singular continuous in this case, by an application of the Lyapunov
exponent criterion for the absence of absolutely continuous spectral
components \cite{BGM,Neil}. The result can now be stated as follows.

\begin{theorem}
  The twisted Fibonacci TDS defined by \eqref{eq:twist-Fib} has mixed
  spectrum of singular type. The pure point part of the diffraction
  measure is\/ $\sum_{k\in L^{\circledast}}I(k)\ts\delta^{}_{k}$, with
  the intensity function\/ $I(k)$ from \eqref{eq:twist-intens}. The
  singular continuous part vanishes if\/
  $u^{}_{a} = u^{}_{\underline{a}}$ and\/
  $u^{}_{b} = u^{}_{\underline{b}}$. \qed
\end{theorem}

Let us now turn to planar systems, which will permit further variations
of geometric origin.

\section{Intermezzo: A planar direct product
  tiling}\label{sec:second}

Here, we start with the obvious in taking the direct product of two
one-dimensional inflation systems to define one for the plane. It is
clear that the spectrum of the resulting system will show the
corresponding direct product structure.

\begin{figure}
\centerline{\includegraphics[width=\textwidth]{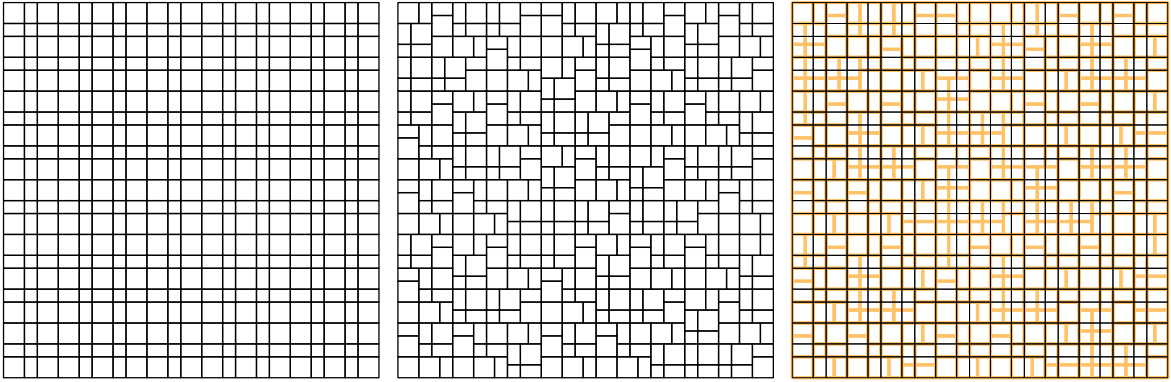}}
\caption{Patches (part of the upper right quadrant of fixed point
  tilings) obtained by six inflation steps from a single large square
  ($3$) by the inflation rule \eqref{eq:fibosqrule} (left) and by
  the modified rule \eqref{eq:scrambled} (centre). The right
  panel shows the superposition of both, with the tiling from the
  modified rule in orange.\label{fig:overlay}}
\end{figure}

For the case of two Fibonacci chains, one thus creates a simple
inflation tiling of the plane with four prototiles, as shown in the
left panel of Figure~\ref{fig:overlay}.  This tiling can be generated
by the inflation rule
\begin{equation}\label{eq:fibosqrule}
  \raisebox{-14pt}{\includegraphics[width=0.8743\textwidth]
    {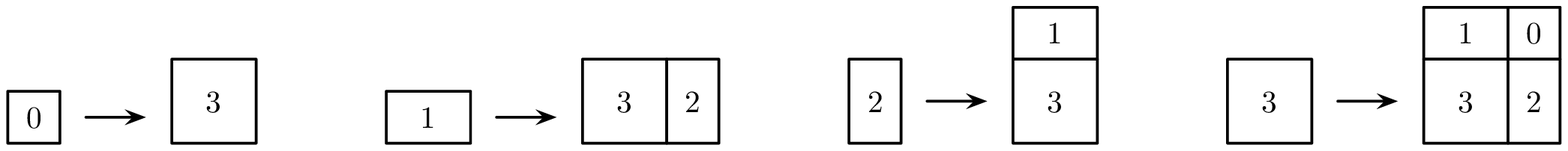}}
\end{equation}
Via iterating from a seed that consists of $4$ big squares, which is
legal, one obtains a sequence that converges towards a $2$-cycle of
infinite planar tilings, either of which can be used to define the
hull $\YY_{2}$ as an orbit closure under the translation action of
$\RR^2$.  The two members of the $2$-cycle are locally
indistinguishable and globally agree on the positive quadrant.

The substitution matrix now reads
\begin{equation}\label{eq:2subst}
    M \, = \, \begin{pmatrix} 
    0 & 0 & 0 & 1 \\
    0 & 0 & 1 & 1 \\
    0 & 1 & 0 & 1 \\
    1 & 1 & 1 & 1 
    \end{pmatrix}  \, = \, 
    \begin{pmatrix} 0 & 1 \\ 1 & 1 \end{pmatrix} \otimes
    \begin{pmatrix} 0 & 1 \\ 1 & 1 \end{pmatrix}
\end{equation}
with PF eigenvalue $\tau^2$ and corresponding eigenvectors
\[
  \langle u\ts | \, = \, \frac{\tau^2}{5}
  \bigl(1,\tau,\tau,\tau^2\bigr)
    \quad\text{and}\quad
    |\ts v\rangle \, = \, \bigl(\tau^{-4},\tau^{-3},
       \tau^{-3},\tau^{-2}\bigr)^{T},
\]
as above normalised such that
$\langle 1|\ts v\rangle = \langle u\ts |\ts v\rangle = 1 $. One
gets
\[
   \lim_{n\to\infty} \tau^{-2n} M^n \, = \, |\ts v\rangle\langle u\ts | 
   \, = \, P \, = \, P^2 ,
\]
where $P$ is again a symmetric projector of rank $1$. 

Any tiling of the hull can be viewed as a Delone set by taking the
lower left corners of the tiles as \emph{control} (or marker) points.
In the positive quadrant, the control point sets of any member of the
above $2$-cycle satisfy the self-similarity relations
\begin{equation}\label{eq:ssr}
\begin{split}
  \vL^{}_{0} &\,=\, \tau \vL^{}_{3} + \twovec{\tau}{\tau}\ts ,\\
  \vL^{}_{1} &\,=\, \tau \vL^{}_{2} + \twovec{0}{\tau} \,\dotcup\,
                   \tau \vL^{}_{3} + \twovec{0}{\tau}\ts ,\\
  \vL^{}_{2} &\,=\, \tau \vL^{}_{1} + \twovec{\tau}{0} \,\dotcup\,
                   \tau \vL^{}_{3} + \twovec{\tau}{0}\ts ,\\
  \vL^{}_{3} &\,=\, \tau \vL^{}_{0} \,\dotcup\, \tau \vL^{}_{1} \,\dotcup\,
                   \tau \vL^{}_{2} \,\dotcup\, \tau \vL^{}_{3}\ts .
\end{split}
\end{equation}
The corresponding displacement matrix is
\[
    T \, = \, \begin{pmatrix} 
     \varnothing & \varnothing & \varnothing & 
     \bigl\{\twovec{\tau}{\tau}\bigr\} \\
     \varnothing & \varnothing & \bigl\{\twovec{0}{\tau}\bigr\} & 
     \bigl\{\twovec{0}{\tau}\bigr\} \\
     \varnothing & \bigl\{\twovec{\tau}{0}\bigr\} & \varnothing & 
     \bigl\{\twovec{\tau}{0}\bigr\} \\
     \bigl\{\twovec{0}{0}\bigr\} & \bigl\{\twovec{0}{0}\bigr\} & 
     \bigl\{\twovec{0}{0}\bigr\} & \bigl\{\twovec{0}{0}\bigr\} 
     \end{pmatrix},
\]
with which the relations \eqref{eq:ssr} simply become
$\vL_{i}=\bigcup_{j}\bigl(\tau \vL_{j} + T_{ij}\bigr)$.

By applying the $\star$-map and taking closures, this turns into
\begin{equation}\label{eq:windowIFS}
  \begin{split}
    W^{}_{0} &\,=\, \sigma W^{}_{3} + \twovec{\sigma}{\sigma}\ts ,\\
  W^{}_{1} &\,=\, \sigma W^{}_{2} + \twovec{0}{\sigma} \,\cup\,
                 \sigma W^{}_{3} + \twovec{0}{\sigma}\ts ,\\
  W^{}_{2} &\,=\, \sigma W^{}_{1} + \twovec{\sigma}{0} \,\cup\,
                 \sigma W^{}_{3} + \twovec{\sigma}{0}\ts ,\\
  W^{}_{3} &\,=\, \sigma W^{}_{0} \,\cup\, \sigma W^{}_{1} \,\cup\,
                 \sigma W^{}_{2} \,\cup\, \sigma W^{}_{3}\ts ,
\end{split}
\end{equation}
where $W^{}_{j}=\overline{\vL^{\star}_{j}}$. Due to taking closures,
the unions need no longer be disjoint. Eq.~\eqref{eq:windowIFS}
defines a contractive IFS on $(\cK \RR^2)^4$, equipped with the
Hausdorff metric topology. It is easy to verify via an explicit
computation that the unique solution --- as expected --- is given by
\begin{equation}\label{eq:windows}
\begin{aligned}
  W^{}_{0} &= [-1,\tau-2]^{2} , &
  W^{}_{1} &= [\tau-2,\tau-1] \times [-1,\tau-2]\ts , \\
  W^{}_{3} &= [\tau-2,\tau-1]^{2} ,
  & W^{}_{2} &= [-1,\tau-2] \times [\tau-2,\tau-1]\ts . 
\end{aligned}
\end{equation}
The $W^{}_{j}$ can now be interpreted as the windows (or rather the
closure of the windows) of the description of the fixed points as
particular projection sets, namely as regular model sets. For the
direct product structure considered here, they are given as products
of the corresponding windows for the two letters in the projection
description of the Fibonacci chain; see the left panel of
Figure~\ref{fig:normal} below for an illustration.

Consequently, the FB coefficients or amplitudes of the defining Delone
sets (for either choice from the $2$-cycle) have product form and are
given by
\begin{align*}
  A^{}_{0} (k) \,& = \, A^{}_{\!\vL_b}\nts
                   (k^{}_{1}) \ts A^{}_{\!\vL_b} (k^{}_{2})
     \, , & 
            A^{}_{1} (k) \, & = \, A^{}_{\!\vL_a}\nts
                              (k^{}_{1}) \ts A^{}_{\!\vL_b} (k^{}_{2})
     \, , \\[1mm]
  A^{}_{2} (k) \, & = \, A^{}_{\!\vL_b}
                    (k^{}_{1}) \ts A^{}_{\!\vL_a}\nts (k^{}_{2})
     \, , & 
            A^{}_{3} (k) \, & = \, A^{}_{\!\vL_a}\nts
                              (k^{}_{1}) \ts A^{}_{\!\vL_a}\nts (k^{}_{2})
     \, ,
\end{align*}
in terms of the one-dimensional amplitudes from Eq.~\eqref{eq:1Dampli},
with $k = (k^{}_{1}, k^{}_{2}) \in \RR^2$.

The internal Fourier matrix is
$\sB(y)=\widecheck{\delta^{}_{T^{\star}}}(y)$, which explicitly reads
\[
   \sB(y) \, = \, \begin{pmatrix}
   0 & 0 & 0 & \ee^{2\pi\ii\sigma (y^{}_1+y^{}_2)} \\
   0 & 0 & \ee^{2\pi\ii\sigma y^{}_2} & \ee^{2\pi\ii\sigma y^{}_2} \\
   0 & \ee^{2\pi\ii\sigma y^{}_1}  & 0 & \ee^{2\pi\ii\sigma y^{}_1} \\
   1 & 1 & 1 & 1
\end{pmatrix} \, = \, 
\begin{pmatrix}
0 & \ee^{2\pi\ii\sigma y^{}_{2}}\\
1 & 1  \end{pmatrix} \otimes
\begin{pmatrix}
0 & \ee^{2\pi\ii\sigma y^{}_{1}}\\
1 & 1  \end{pmatrix} .
\]
This defines the internal cocycle via
$\sB^{(n+1)}(y)=\sB(y)\sB^{(n)}(\sigma y)$ for $n\in\NN$, with
$\sB^{(1)}=\sB$ and $\sB^{(n)}(0)=M^n$. In complete analogy to the
one-dimensional case,
\[ 
   C(y) \, = \, \lim_{n\to\infty} \lvert\sigma\rvert^{2n} \sB^{(n)}(y)
\]
exists, with $C(y)=|c(y)\rangle\langle u|$ and
$|c(0)\rangle=|v\rangle$, so that $C(0)=P$. For any
$k\in L^{\circledast}\!\times\! L^{\circledast}$, where the latter is the Fourier
module (and the dynamical spectrum) of our direct product system, the
FB coefficients are related to $C$ by
\[
   A_{i}(k) \, = \, \myfrac{\tau^2}{5}\,c_{i}(k^{\star})\ts ,
\] 
while they vanish everywhere else. As a consequence, the amplitudes
can be computed, to arbitrary precision, from an approximation of
$C(y)$ via the internal cocycle. To check that this results in the
same expressions for the amplitudes as the direct formula from the
simple rectangular windows is left to the interested reader.  The
result can be summarised as follows.

\begin{theorem}\label{thm:sq-Fib}
  The Fibonacci direct product TDS\/ $(\YY_{\nts 2}, \RR^2)$, in its
  geometric realisation as described above, has pure point spectrum,
  both in the diffraction and in the dynamical sense.  The Fourier
  module is\/ $L^{\circledast}\times L^{\circledast}$ and agrees with
  the dynamical spectrum of\/ $(\YY_{\nts 2}, \RR^2)$ in additive
  notation.  \qed
\end{theorem}

\begin{figure}
\centerline{\includegraphics[width=0.7\textwidth]{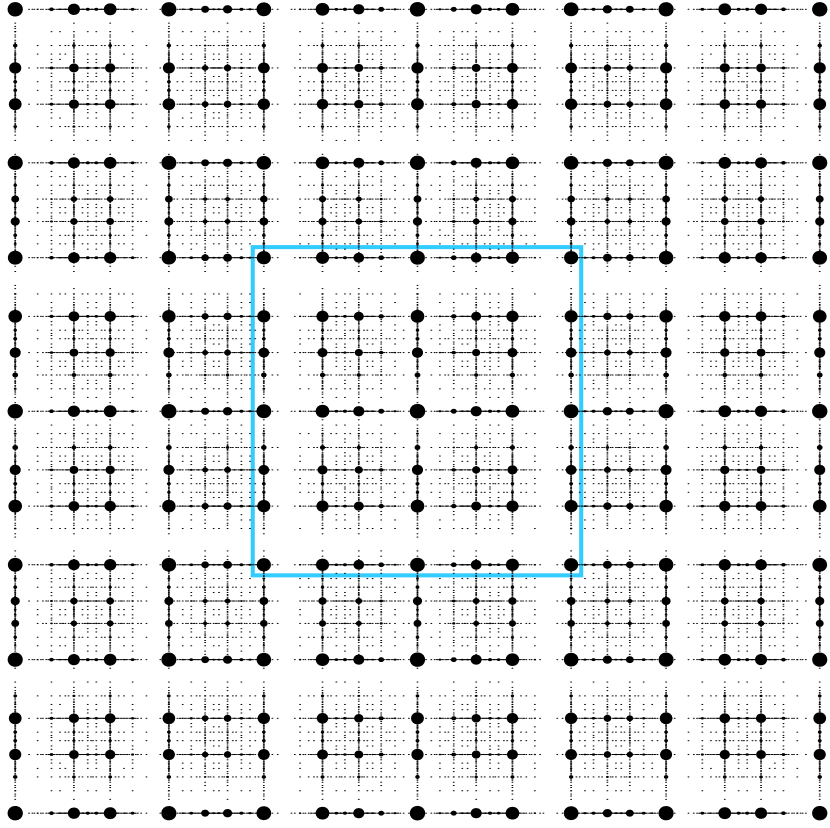}}
\caption{Illustration of the diffraction measure $\widehat{\gamma}$ of
  the Fibonacci direct product system with
  $u^{}_{0}=u^{}_{1}=u^{}_{2}=u^{}_{3}=1$. The individual Dirac
  measures of $\widehat{\gamma}$ for
  $k\in [-5,5]^2 \cap L^{\circledast}\!\times\nts L^{\circledast}$ are
  represented by solid disks centred at the location of the peak, with
  an area that is proportional to the intensity. In later diffraction
  images, we shall only show the central part (blue square) with
  $k\in [-2,2]^2$ for the purpose of convenient
  exposition.\label{fig:squarediff}}
\end{figure}

If we choose weights $u^{}_{0},\dots,u^{}_{3}$ for the four types of
points, the diffraction measure is of the form
$\widehat{\gamma}=\sum_{k\in L^{\circledast}\times L^{\circledast}}
I(k)\,\delta^{}_{k}$ with intensity
$I(k)=\big\lvert\sum_{i=0}^{3} u^{}_{i}\ts
A^{}_{i}(k)\big\rvert^2$. An illustration of $\widehat{\gamma}$ with
all $u^{}_{i}$ equal is shown in Figure~\ref{fig:squarediff}. At this
point, it is a natural task to investigate how the diffraction measure
changes under variations of the direct product structure.

\section{Third variation: Rearranging the direct 
product}\label{sec:third}

Motivated by the DPVs from \cite{Nat-primer,FR}, we now modify the
inflation rule by changing the image of the large square as follows,
\begin{equation}\label{eq:scrambled1}
  \raisebox{-14pt}{\includegraphics[width=0.8743\textwidth]
    {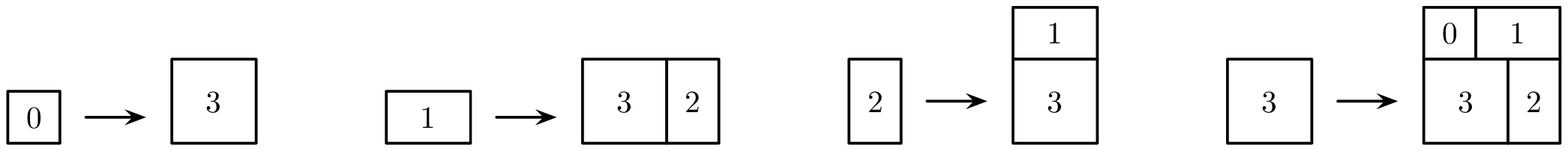}}
\end{equation}
Then, the modified self-similarity relations in the positive quadrant
read
\begin{align*}
  \vL^{\prime}_{0} &\,=\, \tau \vL^{\prime}_{3} + \twovec{0}{\tau}\ts ,\\
  \vL^{\prime}_{1} &\,=\, \tau \vL^{\prime}_{2} +
                     \twovec{0}{\tau} \,\dotcup\,
                   \tau \vL^{\prime}_{3} + \twovec{1}{\tau}\ts ,\\
  \vL^{\prime}_{2} &\,=\, \tau \vL^{\prime}_{1} +
                     \twovec{\tau}{0} \,\dotcup\,
                   \tau \vL^{\prime}_{3} + \twovec{\tau}{0}\ts ,\\
  \vL^{\prime}_{3} &\,=\, \tau \vL^{\prime}_{0}
                     \,\dotcup\, \tau \vL^{\prime}_{1}
                  \,\dotcup\, \tau \vL^{\prime}_{2}
                  \,\dotcup\, \tau \vL^{\prime}_{3}\ts .
\end{align*}
Here, the $\star$-map turns them into the IFS
\begin{equation}\label{eq:windowIFS1}
  \begin{split}
  W^{\prime}_{0} &\,=\, \sigma W^{\prime}_{3} + \twovec{0}{\sigma}\ts ,\\
  W^{\prime}_{1} &\,=\, \sigma W^{\prime}_{2} + \twovec{0}{\sigma} \,\cup\,
                 \sigma W^{\prime}_{3} + \twovec{1}{\sigma}\ts ,\\
  W^{\prime}_{2} &\,=\, \sigma W^{\prime}_{1} + \twovec{\sigma}{0} \,\cup\,
                 \sigma W^{\prime}_{3} + \twovec{\sigma}{0}\ts ,\\
                 W^{\prime}_{3} &\,=\, \sigma W^{\prime}_{0}
                 \,\cup\, \sigma W^{\prime}_{1} \,\cup\,
                 \sigma W^{\prime}_{2} \,\cup\, \sigma W^{\prime}_{3}\ts .
  \end{split}
\end{equation}
It turns out that the closed sets $W^{\prime}_{i}$ that satisfy this IFS
\eqref{eq:windowIFS1} are given by $W^{\prime}_{i}=SW^{}_{i}$ for
$i\in\{0,1,2,3\}$, where
\[
    S \,=\,\begin{pmatrix}1 & -1 \\ 0 & 1\end{pmatrix}.
\]
Indeed, inserting this into the IFS \eqref{eq:windowIFS1} reproduces
the original square Fibonacci IFS \eqref{eq:windowIFS}, except for the
equation for $W^{}_{1}$ which becomes
\[
  W^{}_{1} \,=\, \sigma W^{}_{2} + \twovec{\sigma}{\sigma} \,\cup\,
  \sigma W^{}_{3} + \twovec{1+\sigma}{\sigma}
\]
with different shifts. However, the compact sets $W^{}_{i}$ of
Eq.~\eqref{eq:windows} also satisfy this equation, which,
observing $\sigma < 0$, corresponds
to arranging the rescaled images of $W^{}_{2}$ and $W^{}_{3}$ in the
opposite order, as can be seen by comparing the left panel of
Figure~\ref{fig:normal} with the top left panel of Figure~\ref{fig:quad}.

This relation for the windows translates into a simple relation for
the FB coefficients (or amplitudes). Indeed, for $i\in\{0,1,2,3\}$,
one finds
\begin{equation}\label{eq:ampli1}
  A^{\prime}_{i}(k) \, = \, \widehat{1^{}_{W^{\prime}_{i}}}(-k^{\star})
  \, = \, A^{}_{i}(S^{T}k) \ts ,
\end{equation}
where the second step follows from a simple change of variable
calculation. As can be checked explicitly, $S^{T}$ maps the Fourier
module $L^{\circledast}\!\times\! L^{\circledast}$ onto itself. As a
consequence of this analysis, this DPV also leads to a regular model
set and thus to a system with pure point spectrum, both in the
diffraction and in the dynamical sense, with unchanged
dynamical spectrum.\smallskip
  
Let us now modify the inflation rule by changing the
image of the large square as follows,
\begin{equation}\label{eq:scrambled}
  \raisebox{-14pt}{\includegraphics[width=0.8743\textwidth]
    {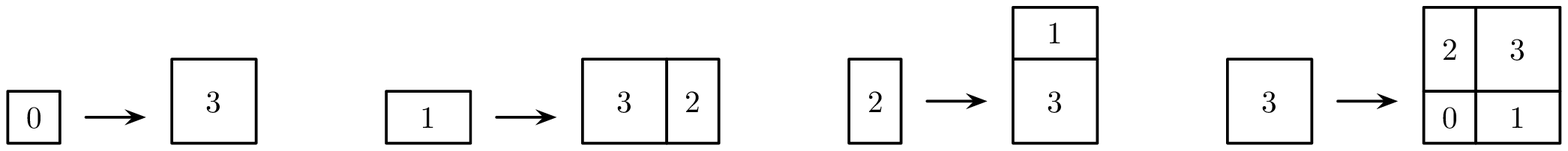}}
\end{equation}
Iterating this inflation rule produces the tiling shown in the central
panel of Figure~\ref{fig:overlay}. Again, the patch consisting of four
large squares is legal (as can be seen in Figure~\ref{fig:overlay}),
and it produces a fixed point under the sixth power of the inflation
rule \eqref{eq:scrambled}. One can determine the corresponding window
IFS in complete analogy to above, but now finds windows with fractal
boundaries; see the top left panel of Figure~\ref{fig:fractwin} below
for an illustration. In fact, these windows are \emph{Rauzy fractals}
\cite{PF,Bernd}, for which one can show that their areas are well
defined and agree with the areas of the windows from
Eq.~\eqref{eq:windows}. Establishing this requires some explicit
estimates of the covering regions on the basis of the contractive IFS,
the details of which we omit here.  Consequently, also this DPV
results in a regular model set.

\begin{prop}\label{prop:castle}
  The dynamical system defined by the DPV rule of
  Eq.~\eqref{eq:scrambled} has pure point spectrum, both in the
  diffraction and the dynamical sense. In particular, the dynamical
  spectrum agrees with that of the Fibonacci direct product system
  from Theorem~\textnormal{\ref{thm:sq-Fib}}, while the
  eigenfunctions, and thus also the diffraction measures, differ.\qed
\end{prop}

\begin{figure}
  \centerline{\includegraphics[width=0.7\textwidth]{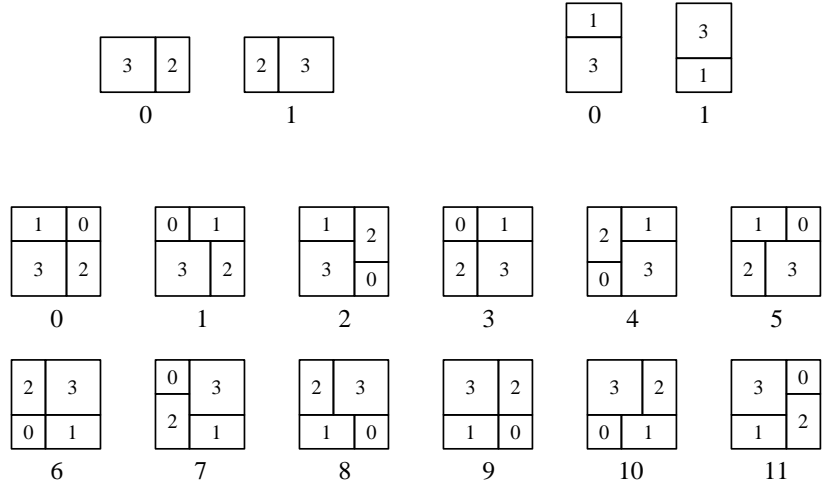}}
  \caption{Labels for the possible decompositions of
    the prototiles of type $1$ and $2$ (top row), and for the $12$
    decompositions of the prototile of type $3$
    (bottom rows).\label{fig:scrambled}}
\end{figure}

\begin{figure}
  \centerline{\includegraphics[width=0.66\textwidth]{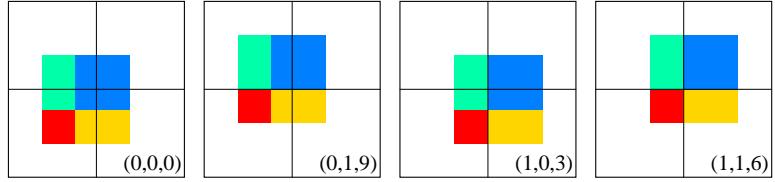}}
  \caption{The four DPVs with the original windows of the square
    Fibonacci tiling, up to a global translation, with parameters as
    shown. The windows for the four types of control points are
    distinguished by colour, namely red ($0$), yellow ($1$), green
    ($2$) and blue ($3$). The outer boxes mark the square
    $[-\tau,\tau]^2$, with the coordinate axes indicated as
    well. Although the fixed points used for the reconstruction of the
    windows are different, these four DPVs lead to the same hull and
    thus define the same dynamical system, namely $(\YY_2, \RR^2)$
    from Theorem~\ref{thm:sq-Fib}.\label{fig:normal}}
\end{figure}

\begin{figure}
  \centerline{\includegraphics[width=0.99\textwidth]{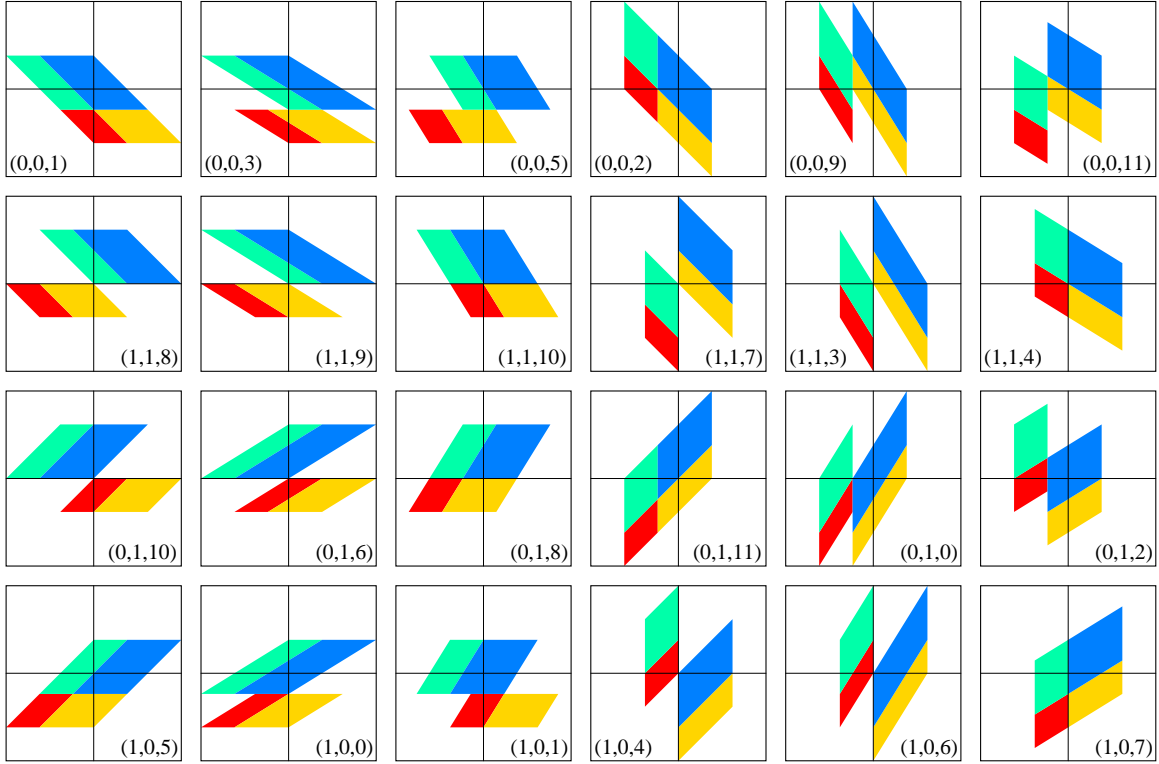}}
  \caption{The $24$ additional cases of DPVs with polygonal windows.
    The left (right) half contains all cases with horizontal
    (vertical) boundaries, grouped columnwise by slope. The windows in
    each column correspond to inflation rules related by reflections,
    with the top and the bottom pair each resulting in tilings that
    are MLD.\label{fig:quad}}
\end{figure}

It is a natural question how different the systems are which are
produced by geometric variations of the Fibonacci direct product.  To
pursue this systematically, we now consider all possible
rearrangements of the stone inflation of each tile type. There is
nothing to rearrange in the inflation of a type-$0$ tile; there are
two rearrangements each of the tiles of type $1$ and $2$, and $12$
possibilities of the type-$3$ tile. These are shown in
Figure~\ref{fig:scrambled}, along with the labelling scheme used.

Note that we only consider rearrangements on the level of the tiles,
and then always use the lower left corner of each prototile as marker
or control point. Other choices are MLD with one of these, and do not
change the spectral type (though they will lead to relatively
translated windows in the cut and project description).  Altogether,
we thus obtain $48$ distinct inflation rules, all with the
substitution matrix from Eq.~\eqref{eq:2subst}. We parameterise the
cases by triples $(i^{}_{1},i^{}_{2},i^{}_{3})$ with
$i^{}_{1},i^{}_{2}\in\{0,1\}$ and $i^{}_{3}\in\{0,1,\ldots,11\}$. In
particular, $(0,0,0)$ is the inflation from Eq.~\eqref{eq:fibosqrule},
while $(0,0,1)$ is the one from Eq.~\eqref{eq:scrambled1} and
$(0,0,6)$ that of Eq.~\eqref{eq:scrambled}.

For each DPV, one can derive the window IFS in the same way as
explained above. For precisely four choices of the parameters, namely
$(0,0,0)$, $(0,1,9)$, $(1,0,3)$ and $(1,1,6)$, one obtains the windows
of Section~\ref{sec:second} or a global translate thereof; see
Figure~\ref{fig:normal}. Observe that the inflation rules $(0,1,9)$
and $(1,0,3)$ emerge from $(0,0,0)$ by a reflection in the horizontal
and vertical axis, respectively, and $(1,1,6)$ by applying both
reflections. Since the original tiling hull is reflection symmetric,
these four rules define the same hull. As a consequence of our
convention to always choose the lower left corner as control point,
which does not preserve the reflection symmetry, it turns out (after
some calculations) that the resulting windows are related by global
shifts. This in particular implies that these four DPVs share the same
diffraction, namely the one illustrated in
Figure~\ref{fig:squarediff}.

Beyond these, there are another $24$ cases where the windows are
parallelograms. They emerge from the original windows by shear
transformations and shifts, as in the example \eqref{eq:scrambled1}
discussed previously. All vertices of the parallelograms are points in
$\ZZ[\tau]^2$ with simple coordinates. These $24$ cases are
illustrated in Figure~\ref{fig:quad}. The proof for each case consists
in determining the window IFS followed by a verification that the
corresponding set of quadrangles satisfies the IFS. We note that
window systems with different slopes cannot produce projection point
set that are MLD, as follows from an application of the general MLD
criterion from \cite[Rem.~7.6]{TAO}. The $24$ window systems thus
partition into $12$ MLD pairs, which are related by double reflections
of the inflation rules, such as $\{(0,0,1),(1,1,8)\}$ or
$\{(0,0,9),(1,1,3)\}$.

\begin{figure}
\centerline{\includegraphics[width=0.9\textwidth]{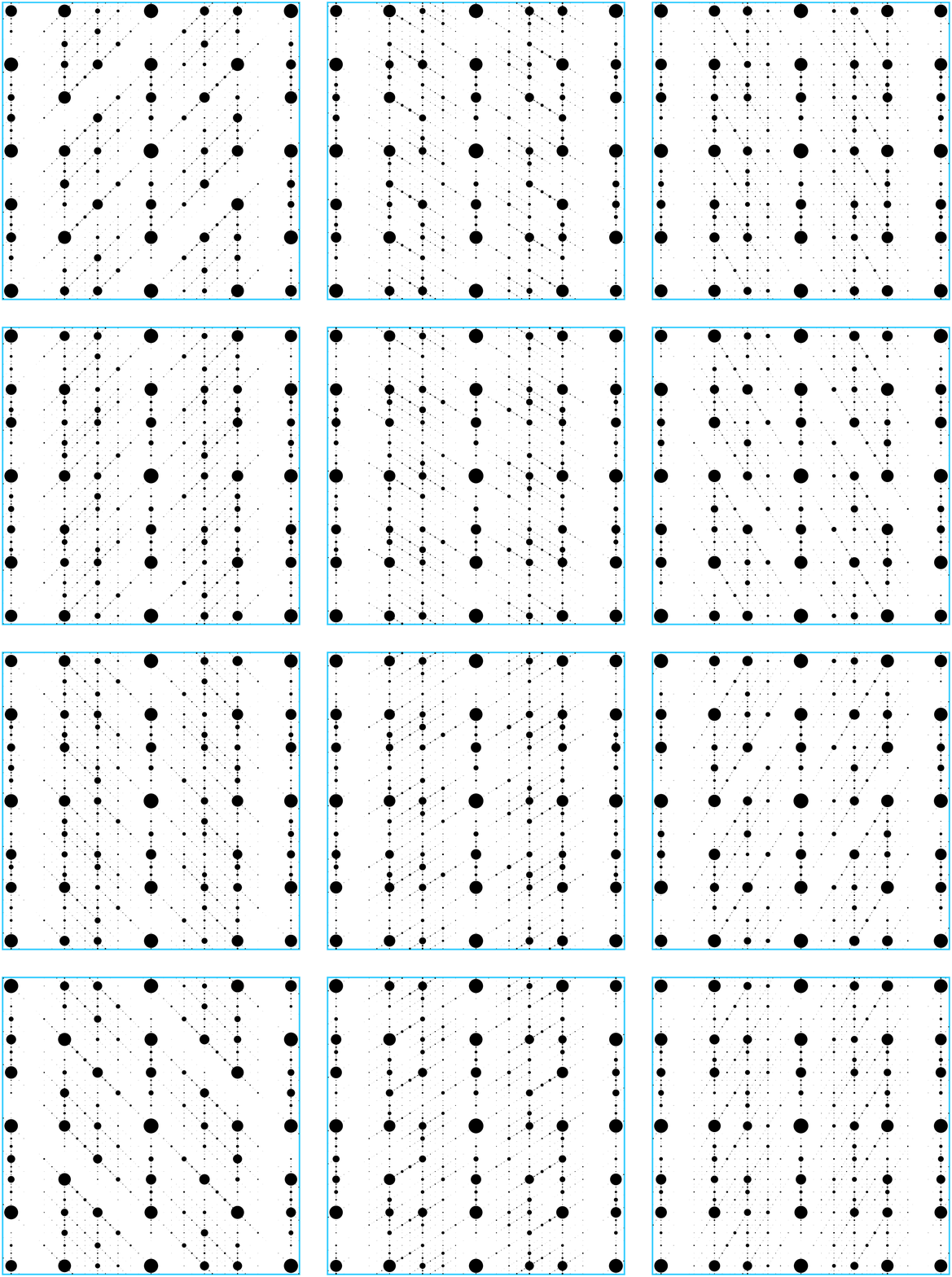}}
\caption{Diffraction measures of the uniform Dirac
  combs for the $12$ DPVs of Figure~\ref{fig:quad} with horizontal
  window boundaries, shown in the same order.\label{fig:quaddiff}}
\end{figure}

\begin{figure}
\centerline{\includegraphics[width=0.85\textwidth]{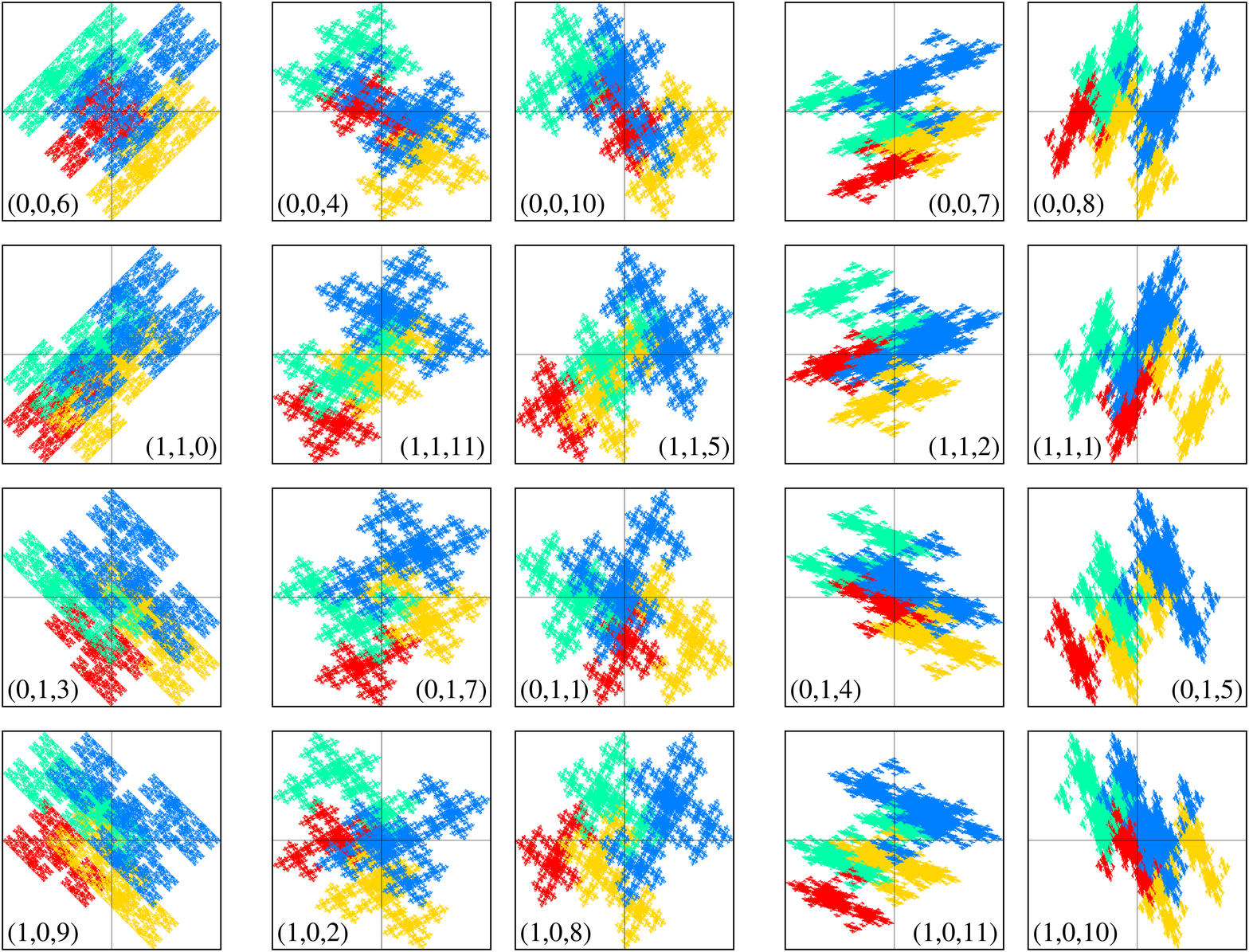}}
\caption{The four Rauzy fractal windows of type `castle' (left
  column), the eight windows of type `cross' (columns 2 and 3) and of
  type `island' (columns 4 and 5). As in Figure~\ref{fig:quad} above,
  the windows in each column correspond to inflation rules that are
  related by reflections. This explains the visible relations in the
  `slopes' of the windows, which are related to certain patterns in
  the diffraction images in
  Figure~\ref{fig:fracdiff}.\label{fig:fractwin}}
\end{figure}

Snapshots of the diffraction measures for the $12$ cases with
horizontal window boundaries are shown in Figure~\ref{fig:quaddiff}.
Let us make some brief comments on their structure, and that of the
other $12$ cases.  Directions in the Fourier module are mapped to
directions in internal space under the $\star$-map, which is totally
discontinuous. The directions of the `white streets' (horizontal or
vertical) are orthogonal to the vertical or horizontal edges of the
windows. The second direction of the window boundaries shows up as a
second direction in the intensity distributions, but in a less obvious
way due to the complicated behaviour of the $\star$-map.\smallskip

The remaining $20$ DPVs lead to windows of Rauzy fractal type. Rauzy
fractals are compact sets that are topologically regular and perfect.
Moreover, they have positive measure and a boundary of fractal (or
Hausdorff) dimension less than that of ambient space, in this case
$d^{}_{\mathrm{H}}<2$; see \cite{Bernd} for a summary of the general
theory. For the Fibonacci DPVs, they come in three topological types,
which we will call `castle', `cross' and `island' to capture their
geometric appearance.

After the preprint of our manuscript became available, Bernd Sing
\cite{Bernie} calculated the Hausdorff dimensions of the three types
of fractals. They all are of the form
\[
     d^{}_{\mathrm{H}} \, = \, \frac{\log (\alpha)}{\log (\tau)}
\]
where $\alpha$ derives from the induced boundary IFS. It is the
largest real root of a fractal-specific integer polynomial, namely
\[
    p (x) \, = \, \begin{cases}
    x^3 - 4 x^2 + 5 x - 3\vspace{1mm} & \text{(castle),} \\
    \begin{array}{@{}l}
        x^{12} - 2 x^{11} - 2 x^{10} + 
       2 x^9 + 4 x^8 - 3 x^7 \\[-1mm]
    \mbox{}\; - 5 x^6 + x^5 + 5 x^4 + 2 x^3 - 2 x^2 - 3 x - 1\vspace{1mm}
    \end{array}
    & \text{(cross),} \\
    x^5 - 2 x^4 - x^3 + 2 x^2 + x - 4 & \text{(island).}
    \end{cases}
\]
The corresponding three dimensions are approximately given 
by  $1.875$, $1.756$ and $1.561$. It is clear from \cite[Rem.~7.6]{TAO}
that projection sets with windows of different Hausdorff dimension
cannot be MLD, which distinguishes the three types from one another,
and also from all cases with polytopal windows.

There are four DPVs with windows of type `castle' (left column of
Figure~\ref{fig:fractwin}). Two of them (the top ones) possess a
reflection symmetry with respect to the diagonal, in the sense that
the windows for the control points of the squares are symmetric, while
the other two are interchanged. The remaining two DPV window systems
form a $2$-cycle under this reflection, with the corresponding
exchange of the control points of the rectangular tiles. These
relations between the DPVs can directly be extracted from the
corresponding inflation rules as well.

The window systems of type `cross' and `island' are also shown in
Figure~\ref{fig:fractwin}. There are eight in each case, whose
interrelations under reflections can be studied on the level of the
inflation rules. Since the windows encode the control points, which
are mapped to new positions under reflections, such reflections result
in more complicated rearrangements of windows and parts of windows. 
An orbit analysis of the variations under rotations and reflections,
as well as the ensuing classification into MLD classes, are left to
the interested reader.

\begin{figure}
\centerline{\includegraphics[width=\textwidth]{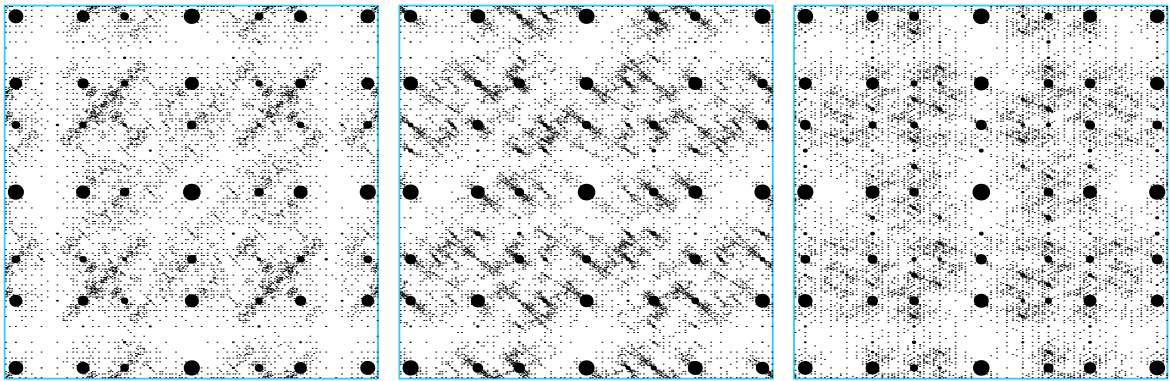}}
\caption{Diffraction measures for selected DPVs with fractal windows,
  $(0,0,6)$ (`castle', left), $(0,0,4)$ (`cross', centre) and
  $(0,0,7)$ (`island', right). While the diffraction image generally
  shows inversion symmetry only, the left image (`castle') displays an
  additional reflection symmetry in the diagonals. This is due to the
  choice of equal weights for the control points of the two
  rectangles.\label{fig:fracdiff}}
\end{figure}

By an inspection analogous to the one that led to
Proposition~\ref{prop:castle}, it is clear that all these examples
lead to regular model sets, some of which have windows with fractal
boundaries. Each individual case comes with its characteristic
internal Fourier matrix, which gives rise to the matrix function
$C(y)=|\ts c(y)\rangle\langle u\ts |$, where $\langle u\ts |$ and
$|\ts c(0)\rangle=|\ts v\rangle$ are always the same, as is the
Fourier module.  Since $A_i(k)=\frac{\tau^{2}}{5} c_{i}(k^{\star})$,
we can compute the FB amplitudes, and hence the diffraction, for all
our examples via the internal cocycle.  Due to the exponentially fast
convergence of the cocycle product, the displayed results are free
from numerical artefacts. Three examples are shown in
Figure~\ref{fig:fracdiff}.

\begin{theorem}
  The\/ $48$ inflation TDSs that emerge from the above DPVs all have
  pure point dynamical spectrum, namely\/ $L^{\circledast} \nts\times
  L^{\circledast}$. These systems are thus measure-theoretically 
  isomorphic by
  the Halmos--von Neumann theorem.  Each individual tiling, via the
  control points, leads to a Dirac comb with pure point diffraction
  measure.\qed
\end{theorem}

It remains an interesting problem to understand the topological
differences between the various DPV classes, which manifest themselves
in different intensity distributions of the diffraction measures, or,
equivalently, in different eigenfunctions for the (commuting) Koopman
operators. In particular, while the $24$ examples of
Figure~\ref{fig:quad} group into $12$ different MLD classes, it is
plausible that they all are topologically conjugate to one another,
and to the Fibonacci direct product TDS. In contrast, it is clear that
the Rauzy fractal windows correspond to new classes under topological
conjugacy, the details of which remain to be analysed.

The above analysis was illustrative, but is of exemplary type
in the sense that it can be applied to any
other Pisot substitution systems. In particular, while we
chose planar examples for ease of presentability, the method works in
higher dimensions as well. In view of the increasingly frequent
appearance of Rauzy fractals as windows, the internal cocycle should
prove to be a valuable tool to gain further insight. Its application is
straightforward when the inflation factor is a unit, but becomes more
delicate when it is not because the internal space is then no longer
Euclidean; see \cite{Bernd} for the general framework required.

\section*{Acknowledgements}

It is our pleasure to thank Lorenzo Sadun, Bernd Sing and Nicolae
Strungaru for helpful discussions.  MB is grateful to Manfred Denker
for his encouragement to analyse geometric and spectral properties of
tilings with methods from dynamical systems theory.  We thank an
anonymous referee for several thoughtful comments, which significantly
helped to improve the presentation. This work was supported by the
German Research Foundation (DFG), within the CRC 1283 at Bielefeld
University, and by EPSRC through grant EP/S010335/1.  


\clearpage


\begin{thebibliography}{111}
\begin{small}

\bibitem{AD}
J.~Aaronson, M.~Denker, O.~Sarig and R.~Zweim\"{u}ller,
Aperiodicity of cocycles and conditional local limit theorems,
\textit{Stoch.\ Dyn.} \textbf{4} (2004) 31--62.

\bibitem{Aki}
S.~Akiyama, M.~Barge, V.~Berth\'{e}, J.-Y.~Lee and A.~Siegel,
On the Pisot substitution conjecture,
in \textit{Mathematics of Aperiodic Order},
eds.\ J.~Kellendonk, D.~Lenz and J.~Savinien
(Birkh\"{a}user, Basel, 2015) pp.~33--72.

\bibitem{BBM}
M.~Baake, M.~Birkner and R.V.~Moody,
Diffraction of stochastic point sets:\ Explicitly
computable examples, \textit{Commun.\ Math.\ Phys.}
\textbf{293} (2010) 611--660;
\texttt{arXiv:0803.1266}.

\bibitem{BFGR}
M.~Baake, N.P.~Frank, U.~Grimm and E.A.~Robinson,
Geometric properties of a binary non-Pisot inflation
and absence of absolutely continuous diffraction,
\textit{Studia Math.} \textbf{247} (2019) 109--154;
\texttt{arXiv:1706.03976}.

\bibitem{BG}
M.~Baake and F.~G\"{a}hler,
Pair correlations of aperiodic inflation rules via 
renormalisation:\ Some interesting examples,
\textit{Topology  \& Appl.} \textbf{205} (2016) 4--27;
\texttt{arXiv:1511.00885}.

\bibitem{BGM}
M.~Baake, F.~G\"{a}hler and N.~Ma\~{n}ibo,
Renormalisation of pair correlation measures for primitive 
inflation rules and absence of absolutely continuous diffraction,
\textit{Commun.\ Math.\ Phys.} \textbf{370} (2019) 591--635;
\texttt{arXiv:1805.09650}.

\bibitem{TAO}
M.~Baake and U.~Grimm,
\textit{Aperiodic Order. Vol.\ 1: A Mathematical Invitation}
(Cambridge University Press, Cambridge, 2013).

\bibitem{BG19}
M.~Baake and U.~Grimm,
Fourier transform of Rauzy fractals and point spectrum
of 1D Pisot inflation tilings,
preprint \texttt{arXiv:1907.11012}.
% update !!

\bibitem{BL}
M.~Baake and D.~Lenz,
Spectral notions of aperiodic order,
\textit{Discr.\ Cont.\ Dynam.\ Syst.\ S}
\textbf{10} (2018) 161--190;
\texttt{arXiv:1601.06629}.

\bibitem{BM-self}
M.~Baake and R.V.~Moody,
Self-similar measures for quasicrystals,
in \textit{Directions in Mathematical Quasicrystals},
eds.\ M.~Baake and R.V.~Moody, CRM Monograph
Series, vol.~13 (AMS, Providence, RI, 2000)
pp.~1--42;
\texttt{arXiv:math.MG/0008063}.

\bibitem{DGS}
M.~Denker, C.~Grillenberger and K.~Sigmund,
\textit{Ergodic Theory of Compact Spaces},
LNM 527 (Springer, Berlin, 1976).

\bibitem{DY}
M.~Denker and M.~Yuri,
Conformal families of measures for general iterated function systems,
in \textit{Recent Trends in Ergodic Theory and Dynamical Systems}, 
eds.\ S.~Bhattacharya, T.~Das, A.~Ghosh and R.~Shah, 
Contemp.\ Math., vol.~631 (AMS, Providence, RI, 2015) 
pp.~93--108.

\bibitem{Nat-primer}
N.P.~Frank,
A primer of substitution tilings of the Euclidean plane,
\textit{Expos.\ Math.} \textbf{26} (2008) 295--326;
\texttt{arXiv:0705.1142}.

\bibitem{FR}
N.P.~Frank and E.A.~Robinson,
Generalized $\beta$-expansions, substitution tilings, and
local finiteness,
\textit{Trans.\ Amer.\ Math.\ Soc.} \textbf{360}  (2008)
1163--1177;
\texttt{arXiv:math.DS/0506098}.

\bibitem{Host}
B.~Host,
Valeurs propres des syst\`{e}mes dynamiques d\'{e}finis par des 
substitutions de longueur variable,
\textit{Ergodic Th.\ \& Dynam.\ Syst.} \textbf{6} (1986) 529--540.

\bibitem{KS}
J.~Kellendonk and L.~Sadun,
Conjugacies of model sets,
\textit{Discr.\ Cont.\ Dynam.\ Syst.\ A}
\textbf{37} (2017) 3805--3830;
\texttt{arXiv:1406.3851}.

\bibitem{Neil}
N.~Ma\~{n}ibo, private communication (2019).

\bibitem{Pohl}
A.D.~Pohl,
Symbolic dynamics, automorphic functions, and Selberg zeta functions
with unitary representations,
in \textit{Dynamics and Numbers}, eds.\ S.~Kolyada, M.~M\"{o}ller,
P.~Moree and T.~Ward, 
Contemp.\ Math., vol.~669 (AMS, Providence, RI, 2016) pp.~205--236;
\texttt{arXiv:1503.00525}.  

\bibitem{PF}
N.~Pytheas Fogg,
\textit{Substitutions in Dynamics, Arithmetics and Combinatorics},
eds.\ V.~Berth\'{e}, S.~Ferenczi, C.~Maduit and A.~Siegel,
LNM 1794 (Springer, Berlin, 2002).

\bibitem{Q}
M.~Queff\'{e}lec,
\textit{Substitution Dynamical Systems --- Spectral Analysis},
2nd ed., LNM 1294 (Springer, Berlin, 2010).

\bibitem{Bernd}
B.~Sing,
\textit{Pisot Substitutions and Beyond},
PhD thesis (Bielefeld University, 2007); 
available electronically at
% \texttt{https://pub.uni-bielefeld.de/record/2302336}.
\texttt{urn:nbn:de:hbz:361-11555}.

\bibitem{Bernie}
B.~Sing, private communication (2019).

\bibitem{Boris}
B.~Solomyak,
Dynamics of self-similar tilings,
\textit{Ergodic Th.\ \& Dynam.\ Syst.} \textbf{17} (1997)
695--738 and \textit{Ergodic Th.\ \& Dynam.\ Syst.}
\textbf{19} (1999) 1685 (Erratum).

\bibitem{Nicu-PC}
N.~Strungaru, private communication (2018).

\bibitem{Weiss}
B.~Weiss,
\textit{Single Orbit Dynamics},
CBMS Regional Conference Series in Mathematics, vol.~95
(AMS, Providence, RI, 2000).

\end{small}
\end{thebibliography}
\end{document}